# CRITERIA FOR IRRATIONALITY OF EULER'S CONSTANT

JONATHAN SONDOW

ABSTRACT. By modifying Beukers' proof of Apéry's theorem that $\zeta(3)$ is irrational, we derive criteria for irrationality of Euler's constant, $\gamma$. For $n > 0$, we define a double integral $I_n$ and a positive integer $S_n$, and prove that with $d_n = \text{LCM}(1,\ldots,n)$ the following are equivalent.
1. The fractional part of $\log S_n$ is given by $\{\log S_n\} = d_{2n} I_n$ for some $n$.
2. The formula holds for all sufficiently large $n$.
3. Euler's constant is a rational number.
A corollary is that if $\{\log S_n\} \geq 2^{-n}$ infinitely often, then $\gamma$ is irrational. Indeed, if the inequality holds for a given $n$ (we present numerical evidence for $1 \leq n \leq 2500$) and $\gamma$ is rational, then its denominator does not divide $d_{2n}\binom{2n}{n}$. We prove a new combinatorial identity in order to show that a certain linear form in logarithms is in fact $\log S_n$. A byproduct is a rapidly converging asymptotic formula for $\gamma$, used by P. Sebah to compute $\gamma$ correct to 18063 decimals.

1. INTRODUCTION

Since R. Apéry's startling 1979 proof [1], [7] that $\zeta(3)$ is irrational, there have been attempts to extend it "up" to $\zeta(n) = \sum_{k=1}^{\infty} k^{-n}$ for odd $n > 3$, with recent progress by K. Ball and T. Rivoal [4] and W. Zudilin [13]. In the present paper, we go "down" to Euler's constant,

$$\gamma := \lim_{N \to \infty} (H_N - \log N), \qquad (1)$$

where $H_N = \sum_{k=1}^{N} k^{-1}$ is the $N^{\text{th}}$ harmonic number, and find criteria for it to be (ir)rational. Note that $\gamma$ can be thought of as "$\zeta(1)$," and as such sits on the fence between rational values of $\zeta(n)$ for $n < 1$, and irrational values (presumably) for $n > 1$.

Defining the double integral

$$I_n := \int_0^1 \int_0^1 -\frac{(x(1-x)y(1-y))^n}{(1-xy)\log xy}\, dx\, dy, \qquad (2)$$

which is inspired by F. Beukers' celebrated integrals for $\zeta(2)$ and $\zeta(3)$ in his elegant 1979 proof [2] of Apéry's theorem, and denoting by $S_n$ the integer product

$$S_n := \prod_{k=1}^{n} \prod_{i=0}^{\min(k-1,n-k)} \prod_{j=i+1}^{n-i} (n+k)^{\frac{2d_{2n}}{j}\binom{n}{i}^2}, \qquad (3)$$

where $d_n = \mathtt{LCM(1,\ldots,n)}$, we prove the following necessary and sufficient conditions for rationality of $\gamma$. (Of course, their negations are then criteria for irrationality of $\gamma$.)

**Rationality Criteria for $\gamma$.** *The following are equivalent*:
(*a*) *The fractional part of* $\log S_n$ *is given by* $\{\log S_n\} = d_{2n} I_n$ *for some* $n$.
(*b*) *The formula holds for all sufficiently large* $n$.
(*c*) *Euler's constant is a rational number.*

We first discovered the Rationality Criteria by modifying the rational function and infinite series in Y. Nesterenko's 1996 proof [6] of Apéry's theorem. In this approach, we defined $I_n$ to be the sum of the series

(4) $$I_n = \sum_{v=n+1}^{\infty} \int_v^{\infty} \left( \frac{n!}{x(x+1)\cdots(x+n)} \right)^2 dx$$

rather than the double integral (2). Remarkably, the two definitions agree. Although we do not know a change of variables transforming one into the other, it turns out that evaluating the integral yields the same result as summing the series, namely,

(5) $$I_n = \binom{2n}{n} \gamma + L_n - A_n$$

where

(6) $$L_n := d_{2n}^{-1} \log S_n \quad \text{and} \quad A_n := \sum_{i=0}^{n} \binom{n}{i}^2 H_{n+i}.$$

(Note that $L_n$ does not actually involve $d_{2n}$, since it cancels out from (3).) We prove (5) for (2) in §2, but defer the proof for (4) to [11], since this series representation of $I_n$ will not be used here. (It might, however, be useful in deciding the arithmetic nature of $\gamma$ from the Rationality Criteria.)

After reading D. Huylebrouck's 2001 survey [5] of multiple integrals in irrationality proofs, we found (2) and rederived the Rationality Criteria from it. It turns out that not only is the integral method more elegant, but also it yields better sufficient conditions for irrationality of gamma – see (*ii*), (*iii*) in Corollary 7 – than the series method does, since the integral is easier to estimate. On the other hand, to prove (5) using the integral definition of $I_n$ requires a combinatorial identity that is not needed if we use the series definition.

An expression like that in (2) gives a new double integral for Euler's constant

$$\gamma = \int_0^1 \int_0^1 -\frac{(1-x)\,dx\,dy}{(1-xy)\log xy}.$$

The proof (see [10]) is similar to that of (5) and is omitted since we will not use the formula.

We prove the Rationality Criteria in §3, and in §4 we give as corollaries several sufficient conditions for irrationality of $\gamma$; they involve $S_n$, but not $I_n$. Here is the most stringent one.

*If $\{\log S_n\} \geq 2^{-n}$ infinitely often, then $\gamma$ is irrational.*

Indeed, the following stronger, finite version holds:

*For fixed $n$, if $\{\log S_n\} \geq 2^{-n}$ and $\gamma$ is a rational number with denominator $q$, then $q$ is not a divisor of $d_{2n}\binom{2n}{n}$; in particular, $|q| > 2n$.*

For example, since $\{\log S_{10000}\} = 0.73...$ (according to P. Sebah [9]) and $0.73 > 2^{-10000}$, if $\gamma$ is rational, its denominator does not divide $d_{20000}\binom{20000}{10000}$.

Table 1 exhibits $S_n$ in the form $S_n = s_n^{2r_n}$, where $r_n := d_n^{-1} d_{2n}$ and $s_n := S_n^{1/(2r_n)}$ $= \exp(\frac{1}{2} d_n L_n)$ are integers; like powers of $n+k$ and $2n-k+1$ are combined.

TABLE 1.

| $n$ | $S_n$ | $\{\log S_n\}$ |
|---|---|---|
| 1 | $2^4$ | 0.7725... |
| 2 | $(12^3)^{12}$ | 0.4566... |
| 3 | $(24^{11} \cdot 5^{38})^{20}$ | 0.3446... |
| 4 | $(40^{25} \cdot 42^{185})^{140}$ | 0.7212... |
| 5 | $(60^{137} \cdot 63^{1762} \cdot 8^{3762})^{84}$ | 0.9645... |
| 6 | $(84^{147} \cdot 88^{2919} \cdot 90^{10794})^{924}$ | 0.5546... |

Computations [9] show that $\{\log S_n\} > 2^{-n}$ for $1 \leq n \leq 2500$. Also, they suggest that $\{\log S_n\}$ is dense in the interval $(0,1)$ and that the cumulative average $<\{\log S(n)\}>$ $= n^{-1} \sum_{k=1}^{n} \{\log S_k\}$ tends to $1/2$ (see Figure 1, courtesy of P. Sebah). To prove $\gamma$ irrational, though, it would suffice just to show that $\{\log S_n\}$ does not tend to zero.

FIGURE 1.

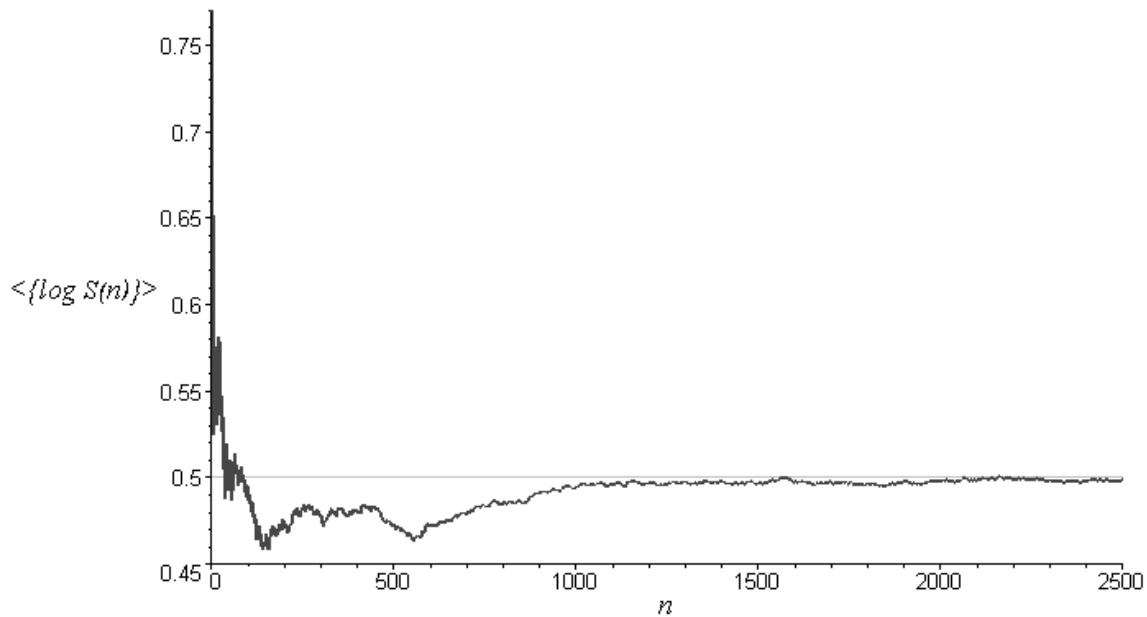

In §5, we prove the asymptotic formula

$$\gamma = \binom{2n}{n}^{-1}(A_n - L_n) + O(2^{-6n} n^{-1/2}).$$

Table 2 shows the first five approximations it provides.

TABLE 2.

| $n$ | $\gamma = 0.577216664901...$ $\binom{2n}{n}^{-1}(A_n - L_n)$ |
|---|---|
| 1 | 0.5(56852819440...) |
| 2 | 0.57(6991219550...) |
| 3 | 0.57721(2786561...) |
| 4 | 0.5772166(25800...) |
| 5 | 0.577216664(353...) |

Using the formula with $n = 10000$, P. Sebah [9] has computed $\gamma$ correct to 18063 decimals, which is precisely the accuracy $2^{-60000}/100$ predicted.

Finally, in the Appendix we establish the combinatorial identity

$$\sum_{i=0}^{k-1}\sum_{j=k}^{n}(-1)^{i+j-1}\binom{n}{i}\binom{n}{j}\frac{1}{j-i} = \sum_{i=0}^{k-1}\binom{n}{i}^2 (H_{n-i} - H_i), \quad 1 \le k \le n,$$

where $H_0 := 0$. Then we use it to prove Lemma 2, which implies that $d_{2n}$ times a certain linear form in logarithms (viz. (11)), a product which is a priori the log of a rational number, is in fact the log of an integer, namely, $\log S_n$. (A reader who wishes to skip the Appendix may take (11) as the definition of $L_n$, and substitute $d_{2n} L_n$ for $\log S_n$ in the Criteria and the corollaries, although this makes their statements less elegant.)

In a paper in preparation [12], we show that conditions stronger than (*i*) and (*ii*) in Corollary 7 imply upper bounds on irrationality measures for $\gamma$.

I thank the referee for suggestions on exposition, and Patrick Gallagher for lectures on [6] and discussions of Laplace's method.

## 2. EVALUATION OF THE DOUBLE INTEGRAL

**Theorem 1.** *The following relation is valid for all $n > 0$:*

(7) $$I_n := \int_0^1\int_0^1 -\frac{(x(1-x)y(1-y))^n}{(1-xy)\log xy}\,dxdy = \binom{2n}{n}\gamma + L_n - A_n,$$

*where $L_n$ is the linear form in logarithms*

(8) $$L_n = \sum_{k=1}^{n}\sum_{i=0}^{\min(k-1,n-k)}\sum_{j=i+1}^{n-i}\binom{n}{i}^2 \frac{2}{j}\log(n+k)$$

*and $d_{2n} A_n \in \mathbf{Z}$.*

*Proof.* Fix $n > 0$. In the integrand of $I_n$, we expand $1/(1-xy)$ in a geometric series with remainder

$$I_n = \int_0^1\int_0^1 -\frac{(x(1-x)y(1-y))^n}{\log xy}\left(\sum_{k=1}^N (xy)^{k-1} + \frac{(xy)^N}{1-xy}\right) dx\,dy$$

and set

(9) $$R_N = \int_0^1\int_0^1 -\frac{(x(1-x)y(1-y))^n (xy)^N}{(1-xy)\log xy}\,dx\,dy.$$

In $I_n - R_N$, we substitute

(10) $$-\frac{1}{\log xy} = \int_0^\infty (xy)^t\,dt,$$

use the Binomial Theorem and reverse the order of integration (which can be justified by replacing $\int_0^1$ by $\int_\varepsilon^{1-\varepsilon}$, and $\int_0^\infty$ by $\int_\varepsilon^{1/\varepsilon}$, and letting $\varepsilon$ tend to $0$), obtaining

$$I_n - R_N = \int_0^1\int_0^1\int_0^\infty \sum_{i,j=0}^n (-1)^{i+j}\binom{n}{i}\binom{n}{j}\sum_{k=1}^N x^{n+i+k-1+t}\,y^{n+j+k-1+t}\,dt\,dx\,dy$$

$$= \int_0^\infty \sum_{i,j=0}^n (-1)^{i+j}\binom{n}{i}\binom{n}{j}\sum_{k=1}^N \frac{1}{(n+i+k+t)(n+j+k+t)}\,dt$$

$$= \int_0^\infty \left(2\sum_{0\le i<j\le n} + \sum_{i=j=0}^n\right)(-1)^{i+j}\binom{n}{i}\binom{n}{j}\sum_{k=1}^N \frac{1}{(n+i+k+t)(n+j+k+t)}\,dt.$$

In this last equality, we have used the symmetry of the expression in $i$ and $j$. It is easily seen that for $0 \leq i < j \leq n$ and $N > 0$, after expansion in partial fractions, the inner sum on $k$ telescopes (cancels) to give

$$\frac{1}{j-i}\sum_{k=1}^{N}\left(\frac{1}{n+i+k+t} - \frac{1}{n+j+k+t}\right) = \frac{1}{j-i}\sum_{k=1}^{j-i}\left(\frac{1}{n+i+k+t} - \frac{1}{N+n+i+k+t}\right).$$

Integrating these terms, as well as the squared ones for $i = j$, yields

$$I_n - R_N = \sum_{i=0}^{n}\binom{n}{i}^2 (H_{N+n+i} - H_{n+i}) + 2\sum_{0 \leq i < j \leq n}\frac{(-1)^{i+j}}{j-i}\binom{n}{i}\binom{n}{j}\sum_{k=1}^{j-i}\log\frac{N+n+i+k}{n+i+k}.$$

Since $n + i \leq 2n$, $n + i + k \leq n + j \leq 2n$ and $n$ is fixed, we have

$$H_{N+n+i} = H_N + o(1), \quad \log(N+n+i+k) = \log N + o(1)$$

as $N \to \infty$. Therefore

$$I_n - R_N = \sum_{i=0}^{n}\binom{n}{i}^2 H_N - 2\sum_{0 \leq i < j \leq n}(-1)^{i+j-1}\binom{n}{i}\binom{n}{j}\log N + L_n - A_n + o(1),$$

where

(11) $$L_n = 2\sum_{0 \leq i < j \leq n}\sum_{k=1}^{j-i}\frac{(-1)^{i+j-1}}{j-i}\binom{n}{i}\binom{n}{j}\log(n+i+k)$$

and $A_n$ is given in (6). Hence $d_{2n} A_n \in \mathbf{Z}$.

Using the identities

$$2 \sum_{0 \leq i < j \leq n} (-1)^{i+j-1} \binom{n}{i}\binom{n}{j} = \sum_{i=0}^{n} \binom{n}{i}^2 = \binom{2n}{n}$$

(to prove the first, put $x = -1$ in $((1+x)^n)^2$; for the second, which is a special case of the Chu-Vandermonde identity, compute the coefficient of $x^n$ in two ways), we obtain that as $N \to \infty$

$$I_n - R_N = \binom{2n}{n}(H_N - \log N) + L_n - A_n + o(1).$$

In view of (1), it only remains to show that $R_N = o(1)$, and that expressions (8) and (11) for $L_n$ agree.

For any $n > 0$, the quantity $\frac{(x(1-x)y(1-y))^n}{1-xy}$ in (9) is bounded by 1 on $[0,1]^2$. Using (10), it follows that

$$0 < R_N < \int_0^1\int_0^1 -\frac{(xy)^N}{\log xy} dx dy = \int_0^\infty \int_0^1 \int_0^1 (xy)^{N+t} dx dy dt,$$

so that $0 < R_N < (N+1)^{-1} \to 0$. The following lemma completes the proof of Theorem 1. •

**Lemma 2.** *The formulas for $L_n$ in (6), (8) and (11) all agree.*

The proof is given in the Appendix.

## 3. RATIONALITY CRITERIA FOR $\gamma$

**Lemma 3.** *The inequalities $1 < d_{2n} < 8^n$ and $0 < I_n < 16^{-n}$ hold for all $n > 0$.*

*Proof.* It is known (see [8], Theorem 13) that $\log d_n < 1.03883\, n$ for all $n$. As $2 \times 1.03883 < 3 \log 2$, we have $1 < d_{2n} < 8^n$.

From (2), observe that $I_n > 0$, and that $t(1-t) < 1/4$ for $t \neq 1/2$ implies that $I_{n+1} < I_n / 16$. Since $I_1 = 2\gamma + 2\log 2 - 5/2 < 1/16$, the lemma is proved. ●

**Theorem 4.** *For $n$ fixed, $\{\log S_n\} = d_{2n} I_n$ if and only if $\gamma = p/q$ for some $p, q \in \mathbf{Z}$, where $q \mid d_{2n}\binom{2n}{n}$.*

*Proof.* Multiply (7) by $d_{2n}$ and use (6) to write the result as

$$\log S_n - d_{2n} I_n = d_{2n} A_n - d_{2n}\binom{2n}{n}\gamma,$$

where $d_{2n} A_n \in \mathbf{Z}$. Note that $\log S_n > 0$, and that $d_{2n} I_n \in (0,1)$, by Lemma 3. Since the fractional part of $x > 0$ is the point $\{x\} \in [0,1)$ such that $x - \{x\} \in \mathbf{Z}$, the theorem follows. ●

**Theorem 5 (Rationality Criteria for $\gamma$).** *The following are equivalent*:
(a) *The fractional part of $\log S_n$ is given by $\{\log S_n\} = d_{2n} I_n$ for some $n$.*
(b) *The formula holds for all sufficiently large $n$.*
(c) *Euler's constant is a rational number.*

*Proof.* Since $q \mid d_{2n}$ for all $n \geq q/2$, Theorem 4 implies that $(c) \Rightarrow (b) \Rightarrow (a) \Rightarrow (c)$. ●

## 4. SUFFICIENT CONDITIONS FOR IRRATIONALITY OF $\gamma$

**Corollary 6.** *If $\{\log S_n\} \geq 2^{-n}$ infinitely often, then $\gamma$ is irrational. In fact, for a given $n$ the inequality implies that if $p, q \in \mathbf{Z}$ and $q \mid d_{2n}\binom{2n}{n}$, then $\gamma \neq p/q$.*

*Proof.* Lemma 3 and Theorem 4 imply the second assertion, which implies the first. ●

**Corollary 7.** *Euler's constant is irrational if at least one of the following is true*:

(i) $\qquad \limsup\limits_{n \to \infty} t^n \{\log S_n\} > 0$ *for some positive number* $t < (4/e)^2$,

(ii) $\qquad \{\log S_n\} \neq (e/4)^{2n(1+o(1))}$,

(iii) $\qquad \lim\limits_{n \to \infty} \dfrac{4^{2n} n}{d_{2n}} \{\log S_n\} \neq \dfrac{\pi}{6 \log 2}$ .

*Proof.* Using $d_n = e^{n(1+o(1))}$ (by the Prime Number Theorem), we see that condition (i) $\Rightarrow$ (ii) $\Rightarrow$ (iii). To show that (iii) $\Rightarrow \gamma \notin \mathbf{Q}$, we make the change of variables $u = 2x - 1$, $v = 2y - 1$ in (2), and find that

$$I_n = 4^{-2n} \int_{-1}^{1} \int_{-1}^{1} \frac{(1-u^2)^n (1-v^2)^n}{(4 - f(u,v)) \log(4 f(u,v)^{-1})} \, du \, dv,$$

where $f(u,v) := (1+u)(1+v)$. Since $1 - w^2 = \exp(-w^2 + O(w^4))$ as $w \to 0$, it follows by Laplace's method (see [3], p. 322) that as $n \to \infty$

(12) $$I_n \sim 4^{-2n} \frac{\pi \, n^{-1}}{3 \log 4},$$

which contradicts (iii) if $\gamma \in \mathbf{Q}$, since then $\{\log S_n\} = d_{2n} I_n$ for $n$ large. This completes the proof. ●

**Corollary 8.** *If $\{\log S_n\} \leq 16 \{\log S_{n+1}\}$ infinitely often when $d_{2n} = d_{2n+2}$, then $\gamma$ is irrational. Indeed, if the inequality holds for a given $n$ and $d_{2n} = d_{2n+2}$, then $\gamma$ is not a rational number with denominator $\leq 2n + 2$.*

*Proof.* The first assertion follows from the second, which in turn follows from Theorem 4, since $I_n > 16 I_{n+1}$. ●

## 5. AN ASYMPTOTIC FORMULA FOR $\gamma$

**Theorem 9.** *The following asymptotic formula holds*:

$$\gamma = \frac{A_n - L_n}{\binom{2n}{n}} + O(2^{-6n} n^{-1/2}).$$

*Proof.* This follows from Theorem 1 together with formula (12) and the estimate $\binom{2n}{n}^{-1} = O(4^{-n} n^{1/2})$ from Stirling's formula. ∎

## APPENDIX. A COMBINATORIAL IDENTITY AND PROOF OF LEMMA 2

**Proposition 10.** *For $1 \leq k \leq n$, we have*

$$(13) \qquad \sum_{i=0}^{k-1} \sum_{j=k}^{n} (-1)^{i+j-1} \binom{n}{i}\binom{n}{j} \frac{1}{j-i} = \sum_{i=0}^{k-1} \binom{n}{i}^2 (H_{n-i} - H_i).$$

We require a lemma.

**Lemma 11.** *We have*

$$(14)_n \qquad \sum_{j=k+1}^{n} (-1)^{k+j-1} \binom{n}{j} \frac{1}{j-k} = \binom{n}{k}(H_n - H_k), \quad 0 \leq k < n,$$

*and*

$$(15)_n \qquad \sum_{i=0}^{k-1} (-1)^{k+i-1} \binom{n}{i} \frac{1}{k-i} = \binom{n}{k}(H_n - H_{n-k}), \quad 0 < k \leq n.$$

*Proof.* We prove $(14)_n$ by induction on $n$. Let $L(n,k)$ denote the sum on the left-hand side of $(14)_n$. For $n=1$, we have $k=0$ and $L(1,0) = 1 = H_1$, as required. Now take $n > 1$. For $k = 0$, we use the recursion

(16) $$\binom{n}{j} = \binom{n-1}{j} + \binom{n-1}{j-1}$$

to break $L(n,0)$ into two sums, in the second of which we substitute

(17) $$\binom{n-1}{j-1}\frac{1}{j} = \binom{n}{j}\frac{1}{n},$$

obtaining

$$L(n,0) = L(n-1,0) + \frac{1}{n}\sum_{j=1}^{n}(-1)^{j-1}\binom{n}{j}$$

$$= L(n-1,0) + \frac{1}{n}.$$

Hence by induction $L(n,0) = H_n$, which proves $(14)_n$ in the case $k=0$.

For $k > 0$, we again use (16) to break $L(n,k)$ into two sums, but now in the second one we set $j' = j - 1$ and write

$$L(n,k) = L(n-1,k) + \sum_{j'=k}^{n-1}(-1)^{k+j'}\binom{n-1}{j'}\frac{1}{j'-(k-1)}$$

$$= L(n-1,k) + L(n-1,k-1).$$

Assuming inductively that $(14)_{n-1}$ holds, we have

$$L(n,k) = \binom{n-1}{k}(H_{n-1} - H_k) + \binom{n-1}{k-1}(H_{n-1} - H_{k-1})$$

$$= \binom{n}{k}(H_{n-1} - H_k) + \binom{n-1}{k-1}\frac{1}{k}$$

and an application of (17) with $j$ replaced by $k$ yields $(14)_n$.

To derive $(15)_n$, we start with $(14)_n$, set $j = n - i$, and replace $k$ by $n - k$, and $\binom{n}{n-m}$ by $\binom{n}{m}$, for $m = i$ and $m = k$. The result is $(15)_n$ and the lemma is proved. ●

*Proof of Proposition* 10. We fix $n$ and induct on $k$. Let $L(k)$ and $R(k)$ denote the left- and right-hand sides of (13). According to $(14)_n$ with $k = 0$, we have $L(1) = H_n = R(1)$. Now observe that

$$L(k+1) - L(k) = \sum_{j=k+1}^{n} (-1)^{k+j-1} \binom{n}{k}\binom{n}{j} \frac{1}{j-k} - \sum_{i=0}^{k-1} (-1)^{i+k-1} \binom{n}{i}\binom{n}{k} \frac{1}{k-i}$$

for $1 \leq k \leq n - 1$. Applying $(14)_n$ and $(15)_n$, we obtain

$$L(k+1) - L(k) = \binom{n}{k}^2 \left(H_{n-k} - H_k\right) = R(k+1) - R(k),$$

which by induction proves (13). ●

We recall Lemma 2 and give the proof.

**Lemma 2.** *The formulas for $L_n$ in* (6), (8) *and* (11) *all agree*.

*Proof.* We may write the right side of (11) as

(18) $$2 \sum_{k=1}^{n} \sum_{i=0}^{k-1} \sum_{j=k}^{n} \binom{n}{i}\binom{n}{j} \frac{(-1)^{i+j-1}}{j-i} \log(n+k)$$

and that of (8) as

(19) $$2 \sum_{k=1}^{n} \sum_{i=0}^{k-1} \binom{n}{i}^2 \left(H_{n-i} - H_i\right) \log(n+k).$$

By Proposition 11, the coefficient of $\log(n+k)$ in (18) is equal to that in (19), for $1 \leq k \leq n$. Hence, the four expressions for $L_n$ in (8), (11), (18) and (19) are all equal. Since those in (6) and (8) evidently agree, we are done. ●

209 West 97th Street, New York, NY 11025
jsondow@alumni.princeton.edu